\documentclass[12pt]{smfart}

\usepackage{amssymb}

\usepackage{smfthm}

\newcommand{\RR}{\mathbf{R}}
\newcommand{\CC}{\mathbf{C}}
\newcommand{\ZZ}{\mathbf{Z}}
\newcommand{\ra}{\rightarrow}
\newcommand{\eps}{\varepsilon}
\newcommand{\capa}{\mathrm{cap}}

\title{L'approximation par des polyn\^omes \`a coefficients entiers}
\date{Avril 2000}

\author{Laurent Berger}
\address{MS 050 Brandeis University \\ Box 549110 \\
Waltham MA 02454-9110}
\email{laurent@brandeis.edu}
\urladdr{http://www.unet.brandeis.edu/\~{}laurent}

\begin{document}

\maketitle
\tableofcontents

\setlength{\parskip}{8pt}
\setlength{\parindent}{0pt}
\setlength{\baselineskip}{15pt}

\section*{Introduction}

Soit $K$ un compact de $\RR$; le th\'eor\`eme de Weierstrass nous dit que
toute fonction $f:K \ra \RR$ continue est limite uniforme d'\'el\'ements de 
$\RR[T]$.

L'objet de cet expos\'e est de d\'eterminer, \'etant donn\'e un compact de
$\RR$, quelles sont les fonctions qui sont limite uniforme d'\'el\'ements de
$\ZZ[T]$. Par exemple, si $0 \in K$ et si une 
telle fonction existe, elle doit \^etre enti\`ere en $0$.

Dans la suite, $K$ d\'esignera un compact de $\RR$ de cardinal infini. Si
$f:K \ra \RR$ est une fonction continue, $|f|_K$ d\'esignera le maximum
de $f$ sur $K$. Un polyn\^ome est dit unitaireunitaire si 
son coefficient dominant vaut $1$.

\section{Compacts de $\RR$ et polyn\^omes de Chebychev}

Commen\c{c}ons par d\'efinir les polyn\^omes de Chebychev d'un compact $K
\subset \RR$.

\begin{theo} 
Soit $K$ un compact de $\RR$ et $n \geq 1$; alors il existe un
polyn\^ome unitaire de degr\'e $n$, not\'e $T_n(K)$, qui r\'ealise
le minimum de $|P_n|_K$ o\`u $P_n$ parcourt l'ensemble des polyn\^omes
unitaires de degr\'e $n$. 
\end{theo}

Ce polyn\^ome s'appelle le $n^{\text{\`eme}}$ polyn\^ome de
Chebychev pour $K$. Si $K=[-1;1]$, on retombe sur les polyn\^omes de
Chebychev classiques (ceci sera d\'emontr\'e plus loin).

\begin{proof}
L'existence vient du fait que dans $\RR_n[T]$, la boule de centre $T^n$ et
de rayon $|T^n|_K$ coupe $\RR_{n-1}[T]$ selon un compact non vide, 
et la fonction $P \mapsto |T^n-P(T)|_K$
y est continue et admet donc un minimum. 
\end{proof}

Le polyn\^ome $T_n(K)$ est unique; pour une d\'emonstration de ce fait, 
voir~\cite[p.140]{tosel}.

\begin{prop}
Si $K=[-1;1]$, alors $T_n(K)=2^{1-n}T_n$, 
les polyn\^omes de Chebychev classiques.
Par suite, 
\[ T_n([a;b])=2 \left( \frac{b-a}{4} \right)^n T_n \left( \frac{2T-a-b}{b-a}
\right) \]
\end{prop}

\begin{proof}
On se ram\`ene \`a la premi\`ere assertion par translation et homoth\'etie.
Rappelons que $T_n$ est d\'efini par $T_n(\cos(\theta))=\cos(n \theta)$.
et que $|T_n|_K$ est r\'ealis\'e par $n+1$ r\'eels de $[-1;1]$. Soit $Q$
de degr\'e $<n$ tel que $|T^n-Q(T)|_K<2^{1-n}$; alors
$2^{1-n}T_n(T)-(T^n-Q(T))$ est un polyn\^ome de degr\'e $<n$ qui s'annule
entre deux extremas cons\'ecutifs de $T_n$ sur $K$, c'est \`a dire en au
moins $n$ points. Il est donc nul. 
\end{proof}

\section{Rayon de capacit\'e des compacts}

Nous allons d\'efinir le rayon de capacit\'e (ou diam\`etre transfini, ou
capacit\'e logarithmique, ou exterior mapping radius) d'un compact.

\begin{prop}
La suite $|T_n(K)|_K^{1/n}$ est convergente; on note $d_1(K)$ sa limite.
\end{prop}

\begin{proof}
Soit $\alpha_n=\log(|T_n(K)|_K^{1/n})$. Si $\alpha_n \ra -\infty$ alors 
$d_1(K)=0$; sinon soit $\alpha= \lim \sup(\alpha_n)$. Comme $T_n(K)T_m(K)$
est un polyn\^ome unitaire de degr\'e $m+n$, on a 
\[ \alpha_{m+n} \leq \alpha_n \frac{n}{n+m} + \alpha_m \frac{m}{n+m} \]
fixons $\eps>0$ et $n$ assez grand. 
On voit que $\alpha_{qn+r} \leq \alpha_n+\eps$
quand $q$ est assez grand ($r$ est entre $0$ et $n$), et donc $\alpha_n \geq
\alpha-\eps$ ce qui montre que la suite $\alpha_n$ converge vers sa limite
sup\'erieure. 
\end{proof}

\begin{prop} 
Soit 
\[ \delta_n(K)=\underset{x_i \in K}{\sup} \prod_{1 \leq i \neq j \leq
n}{|x_i-x_j|^{1/n(n-1)}} \]
alors la suite $\delta_n(K)$ est d\'ecroissante et converge vers un r\'eel
not\'e $d_2(K)$.
\end{prop}

\begin{proof}
On a \[ \delta_{n+1}^{(n-1)n(n+1)}=\prod{|x_i-x_j|^{n-1}}=\prod_k \prod_{1
\leq \hat{k}, i \neq j \leq n}{|x_i-x_j|} \leq \delta_n^{(n-1)n(n+1)}
\] ce qui \'etablit la d\'ecroissance et donc la convergence. 
\end{proof}

\begin{theo}
Les deux constantes $d_1(K)$ et $d_2(K)$ ainsi d\'efinies sont \'egales et
on notera $\capa(K)$ leur valeur commune (rayon de capacit\'e).
\end{theo}

\begin{proof}
Tout d'abord, soient $n$ points $x_i$ qui r\'ealisent le $\sup$ qui d\'efinit
$\delta_n$, et $P(T)=\prod{(T-x_i)}$. On a
\[ \delta_n=\prod{|x_i-x_j|^{1/n(n-1)}}=|\prod P'(x_i)|^{1/n(n-1)} \geq
d_1-\eps \] pour $n$ assez grand ce qui montre que $d_2 \geq d_1$.

Ensuite, on a pour tout $P$ unitaire de degr\'e $n$,  
\[ \delta_{n+1}^{n(n+1)/2} = 
\begin{vmatrix}
1 & \cdots & x_1^{n-1} & P(x_1) \\
\vdots & \ddots & \vdots & \vdots \\
1 & \cdots & x_{n+1}^{n-1} & P(x_{n+1}) 
\end{vmatrix} 
\leq (n+1) \delta_n^{n(n-1)/2} |P|_K \]
comme on le voit en d\'eveloppant le d\'eterminant par rapport \`a la
derni\`ere colonne. Soit $c_n=((n+1)|T_n(K)|_K)^{2/n}$; on trouve
$\delta_{n+1}^{n+1} \leq c_n \delta_n^{n-1}$, et en multipliant ces
in\'egalit\'es pour $n=1, \cdots,k$, on a $\delta_{k+1}^{(k+1)/k}
(\delta_k \cdots
\delta_2)^{1/k} \leq (c_2 \cdots c_k)^{1/k}$. 
On conclut que $d_2 \leq d_1$ en utilisant le th\'eor\`eme 
de Ces\`aro. 
\end{proof}

Par exemple, $\capa([a;b])=(b-a)/4$.

\begin{prop}
Soit $K$ compact; alors $\capa(K) \geq 1$ si et seulement si pour tout
polyn\^ome unitaire $P$ on a $|P|_K \geq 1$. Dans ce cas, $\ZZ[T]$ est 
discret dans $\mathcal{C}^0(K,\RR)$.
\end{prop}

\begin{proof}
S'il existe $P$ unitaire tel que $|P|_K=\alpha<1$ alors $|P^k|_K^{1/k} \leq
\alpha$ et donc $\capa(K)$ aussi.

Soit $f \in \mathcal{C}^0(K,\RR)$ et $P_n$ une suite de polyn\^omes \`a
coefficients entiers qui converge vers $f$. Pour $n>n_0$ assez grand on
aura $|f-P_n|<1/2$, et alors $P_m-P_n$ sera un polyn\^ome entier de norme
$<1$ si $m,n>n_0$, et $P_m-P_n$ divis\'e par son coefficient dominant sera
unitaire de norme $<1$; c'est impossible et donc $P_m=P_n=f$ pour $m,n$ assez
grand. Si de plus $P$ et $Q$ sont distincts \`a coefficients entiers, le
m\^eme argument montre que $|P-Q|_K \geq 1$. 
\end{proof}

\section{Polyn\^omes entiers de petite norme}

On vient de voir que si $\capa(K) \geq 1$, on n'a pas de r\'esultat
int\'eressant d'approximation. \`A partir de maintenant, on va
s'int\'eresser aux compacts $K$ tels que $\capa(K)<1$; la situation est
radicalement diff\'erente. 

Par la proposition pr\'ec\'edente, on dispose d'un polyn\^ome $Q$ unitaire de
norme $<1$.

\begin{prop} Il existe un polyn\^ome $P$ \`a coefficients entiers qui
v\'erifie $|P|_K<1$.
\end{prop}

Cette proposition est vraiment importante, on passe d'une information
analytique ($\capa(K)<1$) \`a une information alg\'ebrique.

\begin{proof}
Soit $\delta>0$, $\alpha=|Q|_K<1$, $d$ le degr\'e de $Q$, $C=1+|T|+\cdots+
|T^{d-1}|$, $\ell_0$ tel que $\alpha^{\ell_0}C/(1-\alpha)<\delta$,
$m=\ell_0d$ et $\eps=\delta/C^{m+1}$.

Soit $k$ assez grand et 
\[ R_k(T)=Q(T)^k-\sum_{\ell \geq \ell_0, i=0 \cdots
d-1}{b_{i,\ell}T^iQ(T)^\ell} \]
o\`u les $b_{i,\ell}$ sont des r\'eels compris entre $0$ et $1$ choisis tels
que l'on puisse \'ecrire $R_k(T)=Z_k(T)+P_k(T)$, avec $Z_k$ \`a coefficients
entiers et $P_k$ de degr\'e $<m$ avec des coefficients entre $0$ et $1$ (un
instant de r\'eflexion montre que c'est toujours possible).

Remarquons que $|R_k-Q^k|_K< \delta$, et que si $k'>k$, $Z_k-Z_{k'}$ est un
polyn\^ome unitaire de degr\'e $k'$ et de norme $|Z_k-Z_{k'}|_K<
|R_k-R_{k'}|_K+|P_k-P_{k'}|_K$. Reste \`a utiliser le principe des tiroirs
pour trouver deux entiers $k$ et $k'$ tels que les coefficients de $P_k$ et
$P_{k'}$ diff\`erent d'au plus $\eps$.

En sommant les erreurs, on trouve que $P=Z_k-Z_{k'}$ est unitaire et entier de
norme $|P|_K<6 \delta$. 
\end{proof}

\section{Noyau de Fekete}

Muni du polyn\^ome $P$ construit pr\'ec\'edemment, nous sommes en mesure
d'approcher des fonctions $f$ v\'erifiant certaines conditions; dans cette
section, nous \'enon\c{c}ons ces conditions. Le compact $K$ est toujours
suppos\'e \^etre de rayon de capacit\'e $<1$. 
On dira que $f:K \ra \RR$ est $\ZZ[T]$-approximable si elle est limite
uniforme sur $K$ de polyn\^omes \`a coefficients entiers.
Si $X \subset K$ est un ensemble, on dit
que $f$ est $X$-interpolable s'il existe un polyn\^ome $R \in \ZZ[T]$ tel
que $f=R$ sur $X$.

Soit $B(K)=\{P \in \ZZ[T], |P|_K<1\}$ (on sait maintenant que $B(K)$ est non
vide), et soit \[ J(K)=\{ x \in K, P(x)=0\ \forall P \in B(K) \}\]
notons que $J(K)$ est fini, car il est contenu dans l'ensemble des z\'eros
d'un polyn\^ome non nul.

\begin{theo}
Soit $K$ un compact tel que $\capa(K)<1$. Alors $f:K \ra \RR$ continue est
$\ZZ[T]$-approximable si et seulement si $f$ est $J(K)$-interpolable.
\end{theo}

\begin{proof}
Si $f$ est $\ZZ[T]$-approximable, alors $P_n \ra f$ et on suppose que 
$|P_n-f|_K<1/2$. Alors $|P_n-P_m|_K<1$, et donc $P_n-P_m$ est nul sur
$J(K)$. Par suite, $f=P_n$ sur $J(K)$.

Pour l'implication contraire, on peut toujours supposer que $f=0$ sur
$J(K)$. Soit $Q_0$ \`a coefficients entiers de norme $<1$. Soient
$x_1,\cdots,x_r$ les z\'eros de $Q_0$ qui sont dans $K$ mais pas dans
$J(K)$~: pour chaque $i$ il existe donc $Q_i \in B(K)$ qui ne s'annule pas
en $x_i$. On pose $Q=\sum_{i \geq 0}Q_i^{2n}$ o\`u $n$ 
est un entier suffisamment
grand. Il est clair que les z\'eros de $Q$ qui sont dans $K$ sont
exactement les \'el\'ements $J(K)$.

On prend $\delta>0$ et $n$ assez grand pour que $\max \{|Q(T)|_K, |TQ(T)|_K
\}<\delta$. Soit $\eps>0$ et $k$ tel que $\sum_{j \geq k}{(j+1)\delta^j}<\eps$.

Soit $K_0$ le compact obtenu en identifiant tous les points de $J(K)$ \`a un
seul, $x_0$. Les fonctions $f$, $Q(T)^k$, et $TQ(T)^k$ sont continues sur
$K_0$. De plus l'alg\`ebre engendr\'ee par $Q(T)^k$ et $TQ(T)^k$ s\'epare
les points de $K_0$. Par le th\'eor\`eme de Stone-Weierstrass, il existe
donc un polyn\^ome \`a deux variables, $\tilde{S}$, tel que
$|f-\tilde{S}(Q(T)^k,TQ(T)^k)|<\eps$, et on peut supposer que le terme constant de
$\tilde{S}$ est nul (car $f(x_0)=0$). Soit $S$ le polyn\^ome obtenu en prenant
les parties enti\`eres des coefficients de $\tilde{S}$. Alors $|S-\tilde{S}|_K<
\sum_{i,j \geq 0,i+j=k}{\delta^{i+j}}<\eps$ et par suite
$S(Q(T)^k,TQ(T)^k)$ approche $f$ \`a $3\eps$ pr\`es. 
\end{proof}

\section{D\'etermination du noyau de Fekete}

Dans cette section, nous indiquons des r\'esultats qui permettent de
simplifier le calcul de $J(K)$; dans la section suivante, nous appliquons
cela au calcul de $J([-a;a])$.

Soit $J_0(K)$ l'ensemble des $\alpha \in J(K)$ qui ont la propri\'et\'e~:
tous les conjugu\'es de $\alpha$ sont r\'eels et appartiennent \`a $K$.

Notre objectif est de d\'emontrer le
\begin{theo}
Les ensembles $J_0(K)$ et $J(K)$ sont \'egaux.
\end{theo}

Pour cela, nous allons d\'emontrer que

\begin{prop}\label{www}
Une fonction continue est $\ZZ[T]$-approximable si et seulement si elle
est $J_0(K)$-interpolable.
\end{prop}

Cela entra\^{\i}ne notamment que $f$ est $J(K)$-interpolable si et seulement si
elle est $J_0(K)$-interpolable, et donc que $J(K)=J_0(K)$.

La preuve de la proposition repose sur le lemme suivant~:

\begin{lemm}
Soit $\{x_1,\cdots,x_r\}$ un ensemble d'entiers alg\'ebriques, 
tel que chacun d'entre eux a un conjugu\'e qui n'est pas dans cet ensemble.
Alors  $\{Q(x_1),\cdots,Q(x_r)\}$, pour $Q$ parcourant $\ZZ[T]$, est dense
dans $\RR^r$.
\end{lemm}

\begin{proof}
On montre tout d'abord le cas o\`u les $x_i$ sont racines d'un m\^eme
polyn\^ome irr\'eductible $P$. Alors soit $x_{r+1}=1$ 
et $V=V(x_i)$ la matrice de Vandermonde construite sur
les $x_i$. Soit $E=\RR^{r+1}$. La matrice $V$ d\'efinit une transformation
lin\'eaire inversible de $E$ dans lui-m\^eme, et l'image 
de $\ZZ^{r+1}$ par $V$ est un r\'eseau de $E$, 
disons $\Lambda$.

Soit $P(R)$ l'ensemble des vecteurs de $E$ dont les $r$ premi\`eres
coordonn\'ees sont de valeur absolue $<1$ et la derni\`ere $<R$. Le
th\'eor\`eme de Minkowski nous fournit, pour $R$ assez grand, un \'el\'ement
non-nul $q \in \Lambda \cap P(R)$. 
Il est facile de voir que $V^{-1}((q_i)_i)$ correspond \`a
un polyn\^ome $Q$ de $\ZZ[T]$ tel que $|Q(x_i)|<1$ pour $i=1 \cdots r$;
enfin $Q(x_i) \neq 0$ pour tout $i$ sinon $Q$ serait nul 
(il est de degr\'e $<$ \`a celui de $P$).

Soient $y_i$ des r\'eels, $k>1$, et $\tilde{P}$ le polyn\^ome de Lagrange 
qui interpole les $y_i/Q(x_i)^k$. Soit $P$ le polyn\^ome dont les
coefficients sont les parties enti\`eres de ceux de $\tilde{P}$.

Alors $|Q^kP(x_i)-y_i| \leq |Q^k(x_i)|(|P(x_i)-y_i/Q(x_i)^k|+
|P(x_i)-\tilde{P}(x_i)|) \leq |Q(x_i)|^kC$ o\`u $C$ ne d\'epend pas de $k$.
Cela \'etablit le r\'esultat (on prend $k$ assez grand).

Si les $x_i$ proviennent de diff\'erents polyn\^omes, alors 
on pose $x_{i,j}$ provenant de $P_j$ irr\'eductible. On se donne $y_{i,j}$
des r\'eels et $\eps>0$.
Soit $Q'_j=\prod_{i \neq j} P_i$. Il existe $Q''_j$ qui
v\'erifie $|Q''_j(x_{i,j})-y_{i,j}/Q'_j(x_{i,j})|<\eps/|Q'_j(x_{i,j})|$. 
Soit alors
$Q=\sum Q'_j Q''_j$. On a $Q(x_{i,j})=Q'_j Q''_j (x_{i,j})$ qui vaut
$y_{i,j}$ \`a $\eps$ pr\`es. 
\end{proof}

\begin{proof}[D\'emonstration de la proposition \ref{www}]
Soit maintenant $J(K)=J_0(K) \cup \{x_1,\cdots,x_r \}$, $\eps>0$, et $P$ le
produit des polyn\^omes minimaux des \'el\'ements de $J_0(K)$. Soit $f$ une
fonction nulle sur $J_0(K)$.

Par le lemme, il existe $Q$ tel que $|Q(x_i)-f(x_i)/P(x_i)|<\eps/|P|_K$. 
Alors $f-QP$ est \`a $\eps$ d'une fonction $g$, nulle sur $J(K)$. 
Comme $g$ est interpolable, il existe $R$ qui l'approche \`a $\eps$ pr\`es
et $QP+R$ approche $f$ \`a $2\eps$ pr\`es. 
\end{proof}

\section{Exemple~: le cas de $[-a;a]$}

Soit $I_a=[-a;a]$. Alors $\capa(I_a)=a/2$. Si $a \geq 2$, il ne se passe
rien d'int\'eressant. Soit donc $a<2$.

Soit $x \in J_0(I_a)$ et $z \in \CC$ tel que $x=z+z^{-1}$. Le complexe $z$
est un entier alg\'ebrique dont tous les conjugu\'es sont de norme $1$. Par
le th\'eor\`eme de Kronecker, c'est une racine de l'unit\'e. Il existe donc
des entiers $j$ et $k$, premiers entre eux, tels que $x=x_j=2\cos(2\pi j/k)$.
Les conjugu\'es de $x$ sont les $x_j$ pour $j \in (\ZZ/k\ZZ)^*$, et doivent
\^etre dans $I_a$ eux aussi, c'est \`a dire que l'on doit avoir
$x_1<a$ ce qui nous donne $k \leq 2\pi/\arccos(a/2)$. 

On a donc~:
\[ J_0([-a;a]) \subset \underset{1 \leq k \leq 
\frac{2\pi}{\arccos(a/2)}}{\cup}
\{2\cos(\frac{2 \pi j}{k}), (j,k)=1 \} \]

Le lecteur est invit\'e \`a traiter le cas des intervalles $[a;b]$ puis
\`a s'essayer \`a des unions disjointes d'intervalles.

\end{document}